\documentclass[12pt,reqno]{amsart}

\newcommand{\R}{{\mathbb R}}

\newcommand{\e}{\epsilon}

\newcommand{\Om}{\Omega}
\newcommand{\D}{\Delta}
\newcommand{\p}{\partial}

\numberwithin{equation}{section}
\newtheorem{theorem}{Theorem}[section]
\newtheorem{defn}[theorem]{Definition}
\newtheorem{lemma}[theorem]{Lemma}
\newtheorem{remark}[theorem]{Remark}

\newtheorem{coro}[theorem]{Corollary}

\begin{document}

\title[Homogenized dynamics of SPDEs ]
{Homogenized dynamics of stochastic partial differential equations
with dynamical boundary conditions }

\author[W. Wang \& J. Duan   ]
{Wei  Wang  \&  Jinqiao Duan }

\address[W.~Wang ]
{Institute of Applied Mathematics\\ Chinese Academy of Sciences\\
Beijing, 100080, China } \email[W.~Wang]{wangwei@amss.ac.cn}

\address[J.~Duan]
{Department of Applied Mathematics\\
Illinois Institute of Technology\\
Chicago, IL 60616, USA} \email[J.~Duan]{duan@iit.edu}

\date{September 13, 2006 submitted; March 2, 2007 accepted}

\thanks{A part of this work was done while J. Duan was
visiting the Ennio De Giorgi Center of Mathematical Research
(www.crm.sns.it), Pisa, Italy. J. Duan would like to thank
Giuseppe Da Prato and Franco Flandoli for their financial support
and hospitality. This work was partly supported by the NSF Grants
DMS-0209326 \& OCE-0620539 and the Outstanding Overseas Chinese
Scholars Fund of the Chinese Academy of Sciences. }


\begin{abstract}
A microscopic heterogeneous system under random influence is
considered. The randomness enters the system at physical boundary of
small scale obstacles as well as at the interior of the physical
medium. This system is modeled by a stochastic partial differential
equation defined on a domain perforated with small holes (obstacles
or heterogeneities), together with random dynamical boundary
conditions on the boundaries of these small holes.

A  homogenized macroscopic model for this microscopic
heterogeneous stochastic system is derived.   This homogenized
effective model is a new stochastic partial differential equation
defined on a unified domain without small holes, with static
boundary condition only. In fact, the random dynamical boundary
conditions are homogenized out, but the impact of random forces on
the small holes' boundaries is quantified as an extra stochastic
term in the homogenized stochastic partial differential equation.
Moreover, the validity of the homogenized model is justified by
showing that the solutions of the microscopic model converge to
those of the effective macroscopic model in probability
distribution, as the size of small holes diminishes to zero.

\end{abstract}

\maketitle

\emph{Dedicated to Giuseppe Da Prato  on the occasion of his 70th
 birthday}

\newpage


\section{Introduction}\label{s1}

Stochastic effects in the multiscale modeling of complex phenomena
have drawn more and more attention in many areas such as material
science \cite{CK97}, climate dynamics \cite{Imkeller}, chemistry
and biology \cite{E00, WaymireDuan}. Stochastic partial
differential equations (SPDEs or stochastic PDEs) arise  naturally
as mathematical models for multiscale systems under random
influences. The need to include stochastic effects in mathematical
modeling of some realistic complex behaviors has become widely
recognized in science  and engineering. But implementing this
approach poses some challenges both in mathematical theory and
computation \cite{Roz, PZ92, HuangYan, Imkeller, WaymireDuan,
Rockner}. The addition
 of stochastic  terms to mathematical models has led to interesting
 new mathematical problems
 at the interface of dynamical systems,  partial
 differential equations, scientific computing, and probability theory.

Sometimes, noise affects a complex system not only inside the
physical medium but also at the physical boundary. Such random
boundary conditions arise in the modeling of, for example, the
air-sea interactions on the ocean surface \cite{PeiOor92}, heat
transfer in a solid in contact with a fluid \cite{Langer},
chemical reactor theory \cite{Lap}, and colloid and interface
chemistry \cite{Vold}.  Random boundary conditions may be static
or dynamical. The \emph{static} boundary conditions, such as
Dirichlet or Neumann boundary conditions, do not involve with time
derivatives of the system state variables. On the contrary, the
\emph{dynamical} boundary conditions contain such time
derivatives. Randomness in such boundary conditions are often due
to various fluctuations.

In this paper we consider a microscopic heterogeneous system,
modeled by a  SPDE  with {\it random dynamical boundary
condition}, in a medium which exhibits small-scale spatial
heterogeneities or obstacles. One example of such microscopic
systems of interest is composite materials containing microscopic
   holes (i.e., cavities), under the impact of   random fluctuations in the
domain and   on the surface of the holes \cite{JKO94, Marchenko}.
A motivation for such a model is based on the consideration that
the interaction between the atoms of the different compositions in
a composite material causes the thermal noise when the scale of
the heterogeneity scale is small. A similar consideration appears
also in  a microscopic stochastic lattice model \cite{BBL05} for a
composite material. Here the microscopic structure is perturbed by
random effect   and the complicated interactions on the boundary
of the holes is     dynamically and randomly evolving. The
heterogeneity scale is assumed to be much smaller than the
macroscopic scale, i.e., we assume that the heterogeneities are
evenly distributed. From a mathematical point of view, one can
assume that microscopic heterogeneities (holes) are periodically
placed in the media. This spatial periodicity with small period
can be represented by a small positive parameter $\e$ (i.e., the
period). In fact we work on the spatial domain $D_\e$,
  obtained by
removing $S_\e$, a collection of small holes of size $\e$,
periodically distributed in a fixed domain $D$. When taking
$\e\rightarrow 0 $, the holes inside domain $D$ are smaller and
smaller and their numbers goes to $\infty$. This signifies that
the heterogeneities are finer and finer.

\medskip

In other words, we consider a spatially extended system with state
variable $u_\e$, where stochastic effects are taken into account
both in the model equation and in the boundary conditions, defined
on a
  domain   perforated with small scale holes.
Specifically, we study a class of stochastic partial differential
equations driven by white noise on a perforated domain with random
dynamical boundary conditions:
\begin{eqnarray*}
du_\e(t, x)&=&\Big[\D u_\e(t,x)+ f(t,x,u_\e, \nabla u_\e)\Big]dt+
g_1(t,x)dW_1(t,x)\\&& in
\;\; D_\e\times (0, T),\nonumber\\
\e^2 du_\e(t,x) &=&\Big[-\frac{\p u_\e(t,x)}{\p
\nu_\e}-\e bu_\e(t,x)\Big]dt+\e g_2(t,x)dW_2(t,x)\\
&&on \;\; \p S_\e\times (0, T).
\end{eqnarray*}
This model will be described in more detail in the next section.


 The goal   is to derive a homogenized  effective
 equation, which is   a new stochastic partial differential equation
 (see Theorems \ref{macro}, \ref{t6.1}, \ref{t6.2} and \ref{t6.3}),  for
the above microscopic heterogenous system, by   homogenization
techniques in the sense of \emph{probability}. Homogenization
theory has been   developed for deterministic systems, and
compactness discussion for the solutions $\{u_\epsilon\}_\epsilon$
in some function space is a key step in various homogenization
approaches \cite{CD99}. However, due to the appearance of the
stochastic terms in the above microscopic system considered in
this paper, such compactness result does not hold for this
stochastic system. Fortunately the compactness  in the sense of
probability, that is, the tightness of the distributions for
$\{u_\epsilon\}$, still holds. So one appropriate way is to
homogenize the stochastic system in the sense of probability. It
is shown that the solution $u_\epsilon$ of the microscopic or
heterogeneous system converges to that of the macroscopic or
homogenized system as $\epsilon\downarrow 0$ in probability
distribution. This means that the distribution  of $\{u_\epsilon
\}_\epsilon$ weakly converges, in some appropriate space, to the
distribution of a stochastic process which solves the macroscopic
effective equation.

It is interesting to note that, for the above system with random
dynamical boundary conditions,   the random force on the boundary
of microscopic scale holes leads, in the homogenization limit, to
a random force distributed all over the physical domain $D$, even
when  the model equation itself contains no stochastic influence
in the domain; see Remark \ref{rmk5.2} in \S 5. We could also say
that the impact of small scale random dynamical boundary
conditions is quantified or carried over to the homogenized model
as an extra random forcing. Therefore, the homogenized effective
model is a new stochastic partial differential equation, defined
on a unified domain without holes.

\medskip

In the present paper, the two-scale convergence techniques are
employed in our approach. Two-scale convergence method is an
important method in homogenization theory which is a formal
mathematic procedure for deriving macroscopic models from
microscopic systems.  Two-scale convergence method contains more
information than the usual weak convergence method; see
\cite{All92} or \S \ref{s4}. Moreover by use of the two-scale
convergence, we do not need the extension operator as introduced
in \cite{CD89}.

Partial differential equations (PDEs) with dynamical boundary
conditions have been studied recently in, for example,
\cite{AKM90, Duncan, Es93, Es95, Hin89, Ti04} and reference
therein. The parabolic SPDEs with noise in the static Neumann
boundary conditions have also been   considered in \cite{PZ92,
PZ96, Mas95}. In \cite{ChSch03}, the authors have studied
well-posedness of the SPDEs with random dynamical boundary
conditions. One of the present authors, with collaborators, has
considered \cite{DuanGaoSchm, YangDuan} dynamical issues of SPDEs
with random dynamical boundary conditions.


The homogenization problem for the deterministic systems defined
in   perforated domains or in other heterogeneous media has been
investigated in, for example, \cite{BOFM92, NR1, NR2, SK02,TTM02}
for heat transfer in a composite material,
\cite{BOFM92,CD89,CDMZ91} for the wave propagation in a composite
material and \cite{LiM05, MP99} for the fluid flow in a porous
media. For a systematic  introduction in homogenization in the
deterministic context, see \cite{CD99, JKO94, SP80, Marchenko}. In
\cite{Ti04}, the effective macroscopic dynamics of a deterministic
partial differential equation   with deterministic dynamical
boundary condition on the microscopic heterogeneity boundary is
studied.

Recently there are also   works on   homogenization of partial
differential equations (PDEs) in the random context; see
\cite{KP02, MM86, PP03, JKO94} for PDEs with random coefficients,
and \cite{BM98, ZV93, ZV94,JKO94} for PDEs in randomly perforated
domains.   A basic assumption in these works is the ergodic
hypotheses on the random coefficients, for the passing of the
limit as $\epsilon \to 0$. Note that the microscopic models in
these works are partial differential equations with random
coefficients, so-called random partial differential equations
(random PDEs)\cite{Souganidis, Kushner,PP03,LiM05,KP02,Wright},
instead of stochastic PDEs --- PDEs with noises --- in the present
paper; see also \cite{WangDuan}. Another novelty in the present
paper is that the microscopic system is under the influence of
random dynamical boundary conditions.


\medskip

We first consider the linear system  and then present   results
about   nonlinear systems with special nonlinear terms. This paper
is organized as follows. The problem formulation is stated in \S
\ref{s2}. Section \ref{s3} is devoted to basic properties of the
microscopic heterogeneous system, and some knowledge  to be used
in our approach is introduced in \S \ref{s4}. The homogenized
effective macroscopic model for the linear system  is derived in
\S \ref{s5}. In the last section,   homogenized effective
macroscopic models  are obtained for three types of nonlinear
systems.


\section{Problem formulation }\label{s2}

Let the physical medium $D$ be an open bounded domain in $\R^n$,
$n\geq 2$, with smooth boundary $\p D$, and let $\e>0$ be a small
parameter. Let $Y=[0,l_1)\times [0,l_2)\times\cdots\times [0,l_n)$
be a representative elementary  cell in $\R^n$ and $S$ an open
subset of $Y$ with smooth boundary $\p S$, such that
$\overline{S}\subset Y$.   The elementary cell $Y$ and the small
cavity or hole $S$ inside it are used to model small scale
obstacles or heterogeneities in a physical medium $D$.  Write
$l=(l_1, l_2, \cdots, l_n)$. Define $\e S=\{\e y: y\in S\}$.
Denote by $S_{\e,k}$ the translated image of $\e S$ by $kl$, $k\in
Z^n$, $kl=(k_1l_1, k_2l_2,\cdots, k_nl_n)$. And let $S_\e$ be the
set all the holes contained in $D$ and $D_\e=D\backslash S_\e$.
Then $D_\e$ is a periodically perforated domain with holes of the
same size as period $\e$. We remark that the holes are assumed to
have no intersection with the boundary $\p D$, which implies that
$\p D_\e=\p D \cup \p S_\e$. See Fig. 1 for the case $n=2$. This
assumption is only needed to avoid technicalities and the results
of
 our paper will remain valid without this assumption  \cite{AMN93}.

\vskip 1cm
\begin{center}
 \setlength{\unitlength}{.35cm}
\begin{picture}(35,13)


  \qbezier(0,2)(0, 1)(1,1)
  \qbezier(0,10)(0,1)(0,2)
  \qbezier(0,10)(0,11)(1,11)
 \qbezier(1,11)(15,11)(15,11)
 \qbezier(15,11)(16,11)(16,10)
 \qbezier(16,10)(16,2)(16,2)
 \qbezier(16,2)(16,1)(15,1)
 \qbezier(15,1)(1,1)(1,1)


\qbezier(.5,1.8)(.5,0.8)(2.5,2.3)
 \qbezier(.5,1.8)(.5,2.6)(.5,2.6)
 \qbezier(.8,2.9)(.5,2.9)(2.5,2.9)
 \qbezier(.5,2.6)(.5,2.9)(.8,2.9)
 \qbezier(2.5,2.9)(3.5,2.9)(2.5,2.3)

 \qbezier(3.5,1.8)(3.5,0.8)(5.5,2.3)
 \qbezier(3.5,1.8)(3.5,2.6)(3.5,2.6)
 \qbezier(3.8,2.9)(3.5,2.9)(5.5,2.9)
 \qbezier(3.5,2.6)(3.5,2.9)(3.8,2.9)
 \qbezier(5.5,2.9)(6.5,2.9)(5.5,2.3)

 \qbezier(6.5,1.8)(6.5,0.8)(8.5,2.3)
 \qbezier(6.5,1.8)(6.5,2.6)(6.5,2.6)
 \qbezier(6.8,2.9)(6.5,2.9)(8.5,2.9)
 \qbezier(6.5,2.6)(6.5,2.9)(6.8,2.9)
 \qbezier(8.5,2.9)(9.5,2.9)(8.5,2.3)

 \qbezier(9.5,1.8)(9.5,0.8)(11.5,2.3)
 \qbezier(9.5,1.8)(9.5,2.6)(9.5,2.6)
 \qbezier(9.8,2.9)(9.5,2.9)(11.5,2.9)
 \qbezier(9.5,2.6)(9.5,2.9)(9.8,2.9)
 \qbezier(11.5,2.9)(12.5,2.9)(11.5,2.3)

 \qbezier(12.5,1.8)(12.5,0.8)(14.5,2.3)
 \qbezier(12.5,1.8)(12.5,2.6)(12.5,2.6)
 \qbezier(12.8,2.9)(12.5,2.9)(14.5,2.9)
 \qbezier(12.5,2.6)(12.5,2.9)(12.8,2.9)
 \qbezier(14.5,2.9)(15.5,2.9)(14.5,2.3)


 \qbezier(.5,3.6)(.5,2.6)(2.5,4.1)
 \qbezier(.5,3.6)(.5,4.4)(.5,4.4)
 \qbezier(.8,4.7)(.5,4.7)(2.5,4.7)
 \qbezier(.5,4.4)(.5,4.7)(.8,4.7)
 \qbezier(2.5,4.7)(3.5,4.7)(2.5,4.1)

 \qbezier(3.5,3.6)(3.5,2.6)(5.5,4.1)
 \qbezier(3.5,3.6)(3.5,4.4)(3.5,4.4)
 \qbezier(3.8,4.7)(3.5,4.7)(5.5,4.7)
 \qbezier(3.5,4.4)(3.5,4.7)(3.8,4.7)
 \qbezier(5.5,4.7)(6.5,4.7)(5.5,4.1)

 \qbezier(6.5,3.6)(6.5,2.6)(8.5,4.1)
 \qbezier(6.5,3.6)(6.5,4.4)(6.5,4.4)
 \qbezier(6.8,4.7)(6.5,4.7)(8.5,4.7)
 \qbezier(6.5,4.4)(6.5,4.7)(6.8,4.7)
 \qbezier(8.5,4.7)(9.5,4.7)(8.5,4.1)

 \qbezier(9.5,3.6)(9.5,2.6)(11.5,4.1)
 \qbezier(9.5,3.6)(9.5,4.4)(9.5,4.4)
 \qbezier(9.8,4.7)(9.5,4.7)(11.5,4.7)
 \qbezier(9.5,4.4)(9.5,4.7)(9.8,4.7)
 \qbezier(11.5,4.7)(12.5,4.7)(11.5,4.1)

 \qbezier(12.5,3.6)(12.5,2.6)(14.5,4.1)
 \qbezier(12.5,3.6)(12.5,4.4)(12.5,4.4)
 \qbezier(12.8,4.7)(12.5,4.7)(14.5,4.7)
 \qbezier(12.5,4.4)(12.5,4.7)(12.8,4.7)
 \qbezier(14.5,4.7)(15.5,4.7)(14.5,4.1)

 \qbezier(.5,5.5)(.5,4.5)(2.5,6)
 \qbezier(.5,5.5)(.5,6.3)(.5,6.3)
 \qbezier(.8,6.6)(.5,6.6)(2.5,6.6)
 \qbezier(.5,6.3)(.5,6.6)(.8,6.6)
 \qbezier(2.5,6.6)(3.5,6.6)(2.5,6)

 \qbezier(3.5,5.5)(3.5,4.5)(5.5,6)
 \qbezier(3.5,5.5)(3.5,6.3)(3.5,6.3)
 \qbezier(3.8,6.6)(3.5,6.6)(5.5,6.6)
 \qbezier(3.5,6.3)(3.5,6.6)(3.8,6.6)
 \qbezier(5.5,6.6)(6.5,6.6)(5.5,6)

 \qbezier(6.5,5.5)(6.5,4.5)(8.5,6)
 \qbezier(6.5,5.5)(6.5,6.3)(6.5,6.3)
 \qbezier(6.8,6.6)(6.5,6.6)(8.5,6.6)
 \qbezier(6.5,6.3)(6.5,6.6)(6.8,6.6)
 \qbezier(8.5,6.6)(9.5,6.6)(8.5,6)

 \qbezier(9.5,5.5)(9.5,4.5)(11.5,6)
 \qbezier(9.5,5.5)(9.5,6.3)(9.5,6.3)
 \qbezier(9.8,6.6)(9.5,6.6)(11.5,6.6)
 \qbezier(9.5,6.3)(9.5,6.6)(9.8,6.6)
 \qbezier(11.5,6.6)(12.5,6.6)(11.5,6)

 \qbezier(12.5,5.5)(12.5,4.5)(14.5,6)
 \qbezier(12.5,5.5)(12.5,6.3)(12.5,6.3)
 \qbezier(12.8,6.6)(12.5,6.6)(14.5,6.6)
 \qbezier(12.5,6.3)(12.5,6.6)(12.8,6.6)
 \qbezier(14.5,6.6)(15.5,6.6)(14.5,6)


 \qbezier(.5,7.3)(.5,6.3)(2.5,7.8)
 \qbezier(.5,7.3)(.5,8.1)(.5,8.1)
 \qbezier(.8,8.4)(.5,8.4)(2.5,8.4)
 \qbezier(.5,8.1)(.5,8.4)(.8,8.4)
 \qbezier(2.5,8.4)(3.5,8.4)(2.5,7.8)

 \qbezier(3.5,7.3)(3.5,6.3)(5.5,7.8)
 \qbezier(3.5,7.3)(3.5,8.1)(3.5,8.1)
 \qbezier(3.8,8.4)(3.5,8.4)(5.5,8.4)
 \qbezier(3.5,8.1)(3.5,8.4)(3.8,8.4)
 \qbezier(5.5,8.4)(6.5,8.4)(5.5,7.8)

 \qbezier(6.5,7.3)(6.5,6.3)(8.5,7.8)
 \qbezier(6.5,7.3)(6.5,8.1)(6.5,8.1)
 \qbezier(6.8,8.4)(6.5,8.4)(8.5,8.4)
 \qbezier(6.5,8.1)(6.5,8.4)(6.8,8.4)
 \qbezier(8.5,8.4)(9.5,8.4)(8.5,7.8)

 \qbezier(9.5,7.3)(9.5,6.3)(11.5,7.8)
 \qbezier(9.5,7.3)(9.5,8.1)(9.5,8.1)
 \qbezier(9.8,8.4)(9.5,8.4)(11.5,8.4)
 \qbezier(9.5,8.1)(9.5,8.4)(9.8,8.4)
 \qbezier(11.5,8.4)(12.5,8.4)(11.5,7.8)

 \qbezier(12.5,7.3)(12.5,6.3)(14.5,7.8)
 \qbezier(12.5,7.3)(12.5,8.1)(12.5,8.1)
 \qbezier(12.8,8.4)(12.5,8.4)(14.5,8.4)
 \qbezier(12.5,8.1)(12.5,8.4)(12.8,8.4)
 \qbezier(14.5,8.4)(15.5,8.4)(14.5,7.8)

 \qbezier(.5,9.1)(.5,8.1)(2.5,9.6)
 \qbezier(.5,9.1)(.5,9.9)(.5,9.9)
 \qbezier(.8,10.2)(.5,10.2)(2.5,10.2)
 \qbezier(.5,9.9)(.5,10.2)(.8,10.2)
 \qbezier(2.5,10.2)(3.5,10.2)(2.5,9.6)

 \qbezier(3.5,9.1)(3.5,8.1)(5.5,9.6)
 \qbezier(3.5,9.1)(3.5,9.9)(3.5,9.9)
 \qbezier(3.8,10.2)(3.5,10.2)(5.5,10.2)
 \qbezier(3.5,9.9)(3.5,10.2)(3.8,10.2)
 \qbezier(5.5,10.2)(6.5,10.2)(5.5,9.6)

 \qbezier(6.5,9.1)(6.5,8.1)(8.5,9.6)
 \qbezier(6.5,9.1)(6.5,9.9)(6.5,9.9)
 \qbezier(6.8,10.2)(6.5,10.2)(8.5,10.2)
 \qbezier(6.5,9.9)(6.5,10.2)(6.8,10.2)
 \qbezier(8.5,10.2)(9.5,10.2)(8.5,9.6)

 \qbezier(9.5,9.1)(9.5,8.1)(11.5,9.6)
 \qbezier(9.5,9.1)(9.5,9.9)(9.5,9.9)
 \qbezier(9.8,10.2)(9.5,10.2)(11.5,10.2)
 \qbezier(9.5,9.9)(9.5,10.2)(9.8,10.2)
 \qbezier(11.5,10.2)(12.5,10.2)(11.5,9.6)

 \qbezier(12.5,9.1)(12.5,8.1)(14.5,9.6)
 \qbezier(12.5,9.1)(12.5,9.9)(12.5,9.9)
 \qbezier(12.8,10.2)(12.5,10.2)(14.5,10.2)
 \qbezier(12.5,9.9)(12.5,10.2)(12.8,10.2)
 \qbezier(14.5,10.2)(15.5,10.2)(14.5,9.6)


 \put(23,6){\vector(-1,0){7.5}}
 \qbezier(15,5.5)(19,3.2)(23,5)
 \put(22.5,4.3){\vector(1,1){1}}

 \put(18.5,6.6){$x=\epsilon\; y$}
 \put(18.5,3.2){$y=\frac{x}{\epsilon}$}


 \put(23,4){\vector(1,0){10}}
 \put(24,2){\vector(0,1){8}}


 \put(33,3){$y_1$}
 \put(22.5,9.5){$y_2$}

 \put(10,12){$D_\epsilon=D\backslash S_\epsilon$}
 \put(10,12){\vector(-2,-1){3}}

 \put(29,2.8){$l_1$}
 \put(23,6.5){$l_2$}
 \put(25,5){$S$}
 \put(30,5.5){\vector(-2,-1){1.7}}
 \put(30,5){$Y^*=Y\backslash\overline{S}$}
 \put(24.5,7.5){$Y=[0,l_1)\times[0,l_2)$}
 \put(22.8,2.8){$O$}

 \qbezier(24,7)(29,7)(29,7)
 \qbezier(29,4)(29,7)(29,7)

 \qbezier(24.5,5.3)(24.5,3.3)(28,6)
 \qbezier(24.5,5.3)(24.5,6.3)(24.5,6.3)
 \qbezier(24.8,6.6)(24.5,6.6)(28,6.6)
 \qbezier(24.5,6.3)(24.5,6.6)(24.8,6.6)
 \qbezier(28,6.6)(29,6.6)(28,6)

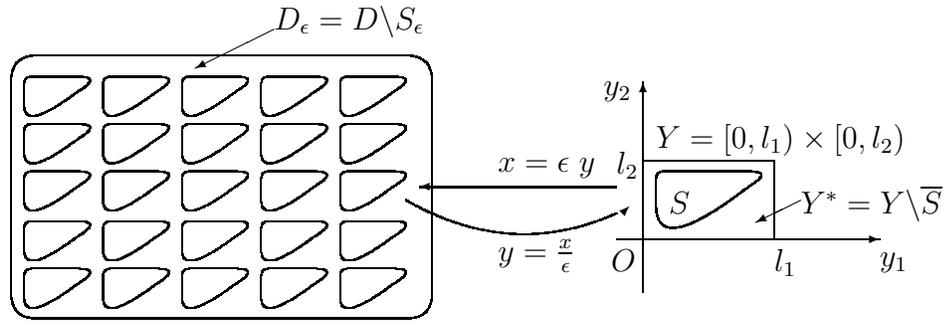
\put(10,-2.5){Fig. 1: Geometric setup in $\mathbb{R}^2$}

\end{picture}
\end{center}

\vskip 1.2cm

In the sequel we use the notations
\begin{equation*}
Y^*=Y\backslash\overline{S}, \;\;\vartheta=\frac{|Y^*|}{|Y|}
\end{equation*}
with $|Y|$ and $|Y^*|$   the  Lebesgue measure of $Y$ and $Y^*$
respectively. Denote by $\chi$ the indicator function, which takes
value 1 on $Y^*$ and value 0 on $Y\setminus Y^*$. In particular,
let $\chi_{A}$ be the indicator function of $A\subset \R^n$. Also
denote by $\tilde{v}$ the zero extension to the whole $D$ for any
function $v$ defined on $D_\e$:
\begin{eqnarray*}
 \tilde{v}= \left\{
  \begin{array}{c l}
     v & \mbox{on  $D_\e$},\\
     0 & \mbox{on $S_\e$}.
  \end{array}
\right.
\end{eqnarray*}
Now for $T>0$ fixed final time,  we consider the following It\^o
type nonautonomous stochastic partial differential equation defined
on the perforated domain $D_\e$ in $\mathbb{R}^n$

\begin{eqnarray}\label{e1}
du_\e(t, x)&=&\Big[\D u_\e(t,x)+ f(t,x,u_\e, \nabla u_\e)\Big]dt+
g_1(t,x)dW_1(t,x)\\&& in
\;\; D_\e\times (0, T),\nonumber\\
\e^2 du_\e(t,x) &=&\Big[-\frac{\p u_\e(t,x)}{\p
\nu_\e}-\e bu_\e(t,x)\Big]dt+\e g_2(t,x)dW_2(t,x)\\
&&on \;\; \p S_\e\times (0,
T),\nonumber \\
u_\e(t,x)&=&0\;\; on\;\; \p D\times (0, T),\\
 u_\e(0,x)&=&u_0(x)\;\; in \;\;
D_\e,\label{e4}
\end{eqnarray}
where $b$ is a real constant, $f: [0, T]\times D\times \R\times
\R^n\rightarrow \R$ satisfies some property which will be
described later and $\nu_\e$ is the exterior unit normal vector on
the boundary $\p S_\e$, $v_0\in L^2(\p S_\e)$ and $u_0\in L^2(D)$.
Moreover, $W_1(t,x)$ and $W_2(t,x)$ are mutually independent
$L^2(D)$ valued
 Wiener processes on a complete probability space $(\Omega,
\mathcal{F}, \mathbb{P})$ with a canonical filtration
$(\mathcal{F}_t)_{t\geq 0}$. Denote by $Q_1$ and $Q_2$ the
covariance operators of $W_1$ and $W_2$ respectively.  Here we
assume that $g_i(t,x)\in\mathcal{L}(L^2(D))$, $i=1,2$ and that
there is a positive constant $C_T$ independent of $\e$ such that
\begin{equation}\label{gi}
\|g_i(t,\cdot)\|_{\mathcal{L}_2^{Q_i}}^2:=\sum_{j=1}^\infty\|g_iQ^\frac{1}{2}_ie_j\|^2_{L^2(D)}\leq
C_T, \;\; i=1,2,\;\;t\in[0, T]
\end{equation}
where $\{e_j\}^\infty_{j=1}$ are eigenvectors  of operator $-\D$
on $D$ with Dirichlet boundary condition and they   form an
orthonormal basis of $L^2(D)$. Here $\mathcal{L}(L^2(D))$ denotes
the space of   bounded linear operators on $L^2(D)$ and
$\mathcal{L}_2^{Q_i}=\mathcal{L}_2^{Q_i}(H)$ denotes the space of
Hilbert-Schmidt operators related to the trace operator $Q_i$
\cite{PZ92}. We also denote by $\mathbf{E}$ the expectation
operator with respect to $\mathbb{P}$.

\medskip

Let $\mathcal{S}$ be a Banach space and $\mathcal{S'}$ be the
strong dual space of $\mathcal{S}$. We recall the definitions and
some properties of weak convergence and $\textrm{weak}^*$
convergence \cite{Yosida}.
\begin{defn}\label{weak1}
A sequence $\{s_n\}$ in $\mathcal{S}$ is said to converge weakly to
$s\in\mathcal{S}$ if $\forall s'\in\mathcal{S'} $,
$$
\lim_{n\rightarrow\infty}(s',s_n)_{\mathcal{S'}, \mathcal{S}}=(s',
s)_{\mathcal{S'}, \mathcal{S}}
$$
which is written as $s_n\rightharpoonup s$ weakly in $\mathcal{S}$.
Note that $(s',s)$ denotes the value of the continuous linear
functional $s'$ at the point $s$.
\end{defn}

\begin{lemma}{\rm (\textbf{Eberlein-Shmulyan})}
Assume that $\mathcal{S}$ is reflexive and let $\{s_n\}$ be a
bounded sequence in $\mathcal{S}$. Then there exists a subsequence
$\{s_{n_k}\}$ and $s\in\mathcal{S}$ such that
$s_{n_k}\rightharpoonup s$ weakly in $\mathcal{S}$ as
$k\rightarrow\infty$. If all the weak convergent subsequence of
$\{s_n\}$ has the same limit $s$, then the whole sequence $\{s_n\}$
weakly converges to $s$.
\end{lemma}

\begin{defn}\label{weak*}
A sequence $\{s'_n\}$ in $\mathcal{S'}$ is said to converge
$weakly^*$ to $s'\in\mathcal{S'}$ if $\forall s\in\mathcal{S} $,
$$
\lim_{n\rightarrow\infty}(s'_n, s)_{\mathcal{S'}, \mathcal{S}}=(s',
s)_{\mathcal{S'}, \mathcal{S}}
$$
which is written as $s'_n\rightharpoonup s'$ $weakly^*$ in
$\mathcal{S'}$.
\end{defn}

\begin{lemma}
Assume that the dual space $\mathcal{S'}$ is reflexive and let
$\{s'_n\}$ be a bounded sequence in $\mathcal{S'}$. Then there
exists a subsequence $\{s'_{n_k}\}$ and $s'\in\mathcal{S'}$ such
that $s'_{n_k}\rightharpoonup s'$ $weakly^*$ in $\mathcal{S'}$ as
$k\rightarrow\infty$. If all the $weakly^*$ convergent subsequence
of $\{s_n'\}$ has the same limit $s'$, then the whole sequence
$\{s'_n\}$ $waekly^*$ converges to $s'$.
\end{lemma}

In the following, for a fixed $T>0$, we always denote by $C_T$ a
constant independent of $\e$. And denote by $D_T$ the set $[0,T]\times D$.  \\


\section{Basic properties of  the microscopic model}\label{s3}

In this section we will present some estimates for solutions of
the microscopic model (\ref{e1}), and then discuss   the tightness
  of the distributions  of the solution processes in some
appropriate space.
We focus our argument in  the case of linear microscopic systems,
where the
  term $f$ is independent of $u_\e$ and $\nabla u_\e$ and
$f(\cdot, \cdot)\in L^2(0, T; L^2(D))$. Then we briefly extend
this to the case of nonlinear microscopic systems with Lipschitz
nonlinearities.

Define by $H^1_\e(D_\e)$  the space of elements of $H^1(D_\e)$
which vanish on $\p D$. Denote by $H_\e^{-1}(D_\e)$ the dual space
of $H_\e^1(D_\e)$ with the usual norm and let $\gamma_\e:
H^1(D_\e) \rightarrow  L^2(\p S_\e)$ be  the trace operator with
respect to $\p S_\e$ which is continuous  \cite{Tri78}. We also
denote that $H^\frac{1}{2}(\p S_\e)= \gamma_\e (H^1(D_\e))$ and
let $H_\e^{-\frac12}(D_\e)$ be the dual space of $H_\e^\frac12
(D_\e)$.

Introduce the following function spaces
$$
X_\e^1=\big\{(u, v)\in H^1_\e(D_\e)\times H_\e^\frac{1}{2}(\p S_\e):
v=\e \gamma_\e u \big\}
$$
and
$$
X_\e^0=\big\{L^2(D_\e)\times L_\e^2(\p S_\e)\big\}
$$
with the usual product and norm. Define an operator $B_\e$ on the
space $H_\e^1(D_\e)$ as
\begin{equation}\label{B1}
B_\e u=\frac{\p u}{\p n}+\e bu,\;\; u\in H_\e^1(D_\e).
\end{equation}
Now we define the operator $\mathcal{A}_\e$ on $\mathcal{D(A_\e)}=\{
(u, v)\in X_\e^1: (-\D u, R_\e B_\e u)\in X_\e^0\}$, where $R_\e$ is
the restriction to $\p S_\e$, as
\begin{equation}\label{Ae}
\mathcal{A}_\e z=(-\D u, \frac{1}{\e}R_\e B_\e u), \;\; z=(u,
v)\in\mathcal{D(A_\e)}.
\end{equation}

Associated with the operator $\mathcal{A}_\e$, we introduce the
bilinear form on $X_\e^1$
\begin{equation}\label{ae}
a_\e(z,\bar{z})=\int_{D_\e} \nabla u\nabla \bar{u}dx+\e
b\int_{S_\e}\gamma_\e(u)\gamma_\e(\bar{u})ds
\end{equation}
with $z=(u, v)$, $\bar{z}=(\bar{u},\bar{v})\in X_\e^1$.  Notice
that $|\gamma_\e(u)|^2_{L^2(\p S_\e)}\leq
C(S_\e)|u|^2_{H_{\e}^1(D_\e)}$, we see that there is $M>0$,
independent of $\e$, such that
$$
a_\e(z,\bar{z})\leq M|u|_{H_\e^1(D_\e)}|\bar{u}|_{H_\e^1(D_\e)}
$$
and the following coercive property of $a$ holds
\begin{equation}\label{coercive}
a_\e(z,z)\geq
\bar{\alpha}|z|^2_{X_\e^1}-\bar{\beta}|z_\e|^2_{X_\e^0},\;\; z\in
X^1_\e
\end{equation}
for some constants $\bar{\alpha}$, $\bar{\beta}>0$ which are also
independent of $\e$. Write the $C_0$-semigroup generated by
operator $-\mathcal{A}_\e$ as $S_\e(t)$.

 Then the system (\ref{e1})-(\ref{e4})
can be rewritten as the following abstract stochastic evolutionary
equation
\begin{equation}\label{abse}
dz_\e(t,x)=[-\mathcal{A}_\e
z_\e(t,x)+F_\e(t,x)]dt+G_\e(t,x)dW(t,x),\;\; z_\e(0)=z_0
\end{equation}
where
$$
F_\e(t,x)=(f(t,x),0)^t,\;\;G_\e(t,x)dW(t)=(g_1(t,x)dW_1(t,x),
g_2(t,x)dW_2(t,x))^t
$$
and $z_0=(u_0, v_0)$. And the solution of (\ref{abse}) can  be
written in the mild sense
\begin{equation}\label{mild}
z_\e(t)=S_\e(t)z_0+\int_0^tS_\e(t-s)F_\e(s)ds+\int_0^tS_\e(t-s)G_\e(s)dW(s).
\end{equation}
  Moreover, the variational formulation is
\begin{eqnarray}\label{variational}
&&\int_0^T\hspace{-.2cm}\int_{D_\e}\dot{u}_\e\varphi
dxdt+\e^2\int_0^T\hspace{-.2cm}\int_{\p S_\e}\dot{u}_\e\varphi
dxdt+\e
b\int_0^T\hspace{-.2cm}\int_{\p S_\e}u_\e \varphi dxdt=\\
&\hspace{-.2cm}-&\hspace{-.2cm}\int_0^T\hspace{-.2cm}\int_{D_\e}\nabla
u_\e\nabla \varphi dxdt+ \int_0^T\hspace{-.2cm}\int_{D_\e}f\varphi
dxdt+\int_0^T\hspace{-.2cm}\int_{D_\e}g_1\varphi\dot{W}_1dxdt+\e\int_0^T\hspace{-.2cm}\int_{\p
S_\e}g_2\varphi\dot{W}_2dxdt \nonumber
\end{eqnarray}
for $\varphi(t, x)\in C_0^\infty([0, T]\times D_\e)$. Here $\dot{}$
denotes $\frac{d}{dt}$.


\medskip

For the well-posedness of system (\ref{abse}) we have the following
result.

\begin{theorem}\label{wellpose} {\rm (\textbf{Global well-posedness of microscopic
model})} Assume that (\ref{gi}) holds for   $T>0$. If $z_0=(u_0,
v_0)$ is  a $\big(\mathcal{F}_0,
\mathcal{B}(X^0_\e)\big)$-measurable random variable, then the
system (\ref{abse}) has a unique mild solution $z_\e\in
L^2\big(\Omega, C(0,T;$ \\ $X^0_\e)\cap L^2(0, T; X^1_\e)\big)$,
which is also a weak solution in the following sense
\begin{eqnarray}
&&\hspace{-0.3cm} (z_\e(t),\phi)_{X^1_\e}\label{weak}\\
\hspace{-0.3cm}&=&\hspace{-0.3cm}(z_0,
\phi)_{X^1_\e}+\int_0^t(-\mathcal{A}_\e z_\e(s),
\phi)_{X^1_\e}ds+\int_0^t(F_\e, \phi)_{X^1_\e}ds+ \int_0^t(G_\e dW,
\phi)_{X^1_\e}\nonumber
\end{eqnarray}
for $t\in[0,T)$ and $\phi\in X^1_\e$. Moreover if $z_0$ is
independent of $W(t)$ with $\mathbf{E}|z_0|^2_{X^0_\e} <\infty $,
then
\begin{equation}\label{est1}
\mathbf{E}|z_\e(t)|^2_{X^0_\e}+\mathbf{E}\int_0^t|z_\e(s)|^2_{X^1_\e}ds\leq
(1+\mathbf{E}|z_0|^2_{X^0_\e})C_T,\;\; {\rm for} \;\; t\in [0, T]
\end{equation}
and
\begin{equation}\label{est2}
\mathbf{E}\big\{\sup_{t\in[0, T]}|z_\e(t)|^2_{X^0_\e}\big\}\leq
\big(1+\mathbf{E}|z_0|^2_{X^0_\e}+\mathbf{E}\int_0^T|z_\e(s)|^2_{X_\e^1}ds\big)C_T.
\end{equation}
\end{theorem}

\begin{proof}
By the assumption (\ref{gi}), we have
$$
 \|G_\e(t, x)\|^2_{\mathcal{L}^Q_2}=\|g_1(t, x)\|^2_{\mathcal{L}^{Q_1}_2}+\|g_2(t,
 x)\|^2_{\mathcal{L}^{Q_2}_2}<\infty.
$$
Then the classical result   \cite{PZ92} yields the local existence
of $z_\e$. By applying the stochastic Fubini theorem \cite{PZ92},
it can be verified that the local mild solution is also a weak
solution.

Now we give the following a $priori$ estimates which yields the
existence of weak solution on $[0 ,T]$ for any $T>0$.

Applying It\^o formula to $|z_\e|_{X_\e^0}^2$, we derive
\begin{eqnarray}\label{It1}
d|z_\e(t)|^2_{X^0_\e} +2(\mathcal{A}_\e z_\e,
z_\e)_{X^0_\e}dt&=&2(F_\e(t,x), z_\e)_{X^0_\e}dt+2(G_\e(t,x)dW(t),
z_\e)_{X^0_\e}+\nonumber\\ &&|G_\e(t,x)|^2_{\mathcal{L}_2^Q}dt.
\end{eqnarray}
By the coercivity (\ref{coercive}) of $a_\e(\cdot,\cdot)$,
 integrating (\ref{It1}) with respect to $t$ yields
\begin{eqnarray*}
&&|z_\e(t)|^2_{X^0_\e} +2\bar{\alpha} \int_0^t|z_\e(s)|^2_{X^1_\e}ds\\
&\leq &|z_0|^2_{X^0_\e}+|F_\e|^2_{L^2(0,T;
X^0_\e)}+(2\bar{\beta}+1)\int_0^t|z_\e(s)|^2_{X_\e^0}ds+\\
&&2\int_0^t(G_\e(s) dW(s),
z_\e(s))_{X^0_\e}+\int_0^t|G_\e(s)|^2_{\mathcal{L}_2^Q}ds.
\end{eqnarray*}
Taking expectation on both sides of the above inequality yields
\begin{eqnarray*}
&&\mathbf{E}|z_\e(t)|^2_{X^0_\e} +2\bar{\alpha} \mathbf{E}\int_0^t|z_\e(s)|^2_{X^1_\e}ds\\
&\leq &\mathbf{E}|z_0|^2_{X^0_\e}+|F_\e|^2_{L^2(0,T;
X^0_\e)}+(2\bar{\beta}+1)\int_0^t\mathbf{E}|z_\e(s)|^2_{X_\e^0}ds+\int_0^t|G_\e(s)|^2_{\mathcal{L}_2^Q}ds.
\end{eqnarray*}
Then the  Gronwall lemma gives the estimate (\ref{est1}). Notice
that, by Lemma 7.2 in \cite{PZ92},
\begin{equation*}
\mathbf{E}\sup_{t\in [0,
T]}\Big|\int_0^tS_\e(t-s)G_\e(s,x)ds\Big|^2_{X_\e^0}\leq
C_T\int_0^T|G_\e(s)|^2_{\mathcal{L}_2^Q}ds.
\end{equation*}
Therefore by the assumption on $f$ and (\ref{mild}) we have the
estimate (\ref{est2}). The proof is hence complete.
\end{proof}

By the above result and the definition of $z_\e$ we have the
following corollary.
\begin{coro}\label{est0}
Assume the conditions in Theorem \ref{wellpose}. Then for $t\in
[0, T]$, we have
\begin{eqnarray}\label{est3}
&& \mathbf{E}\big(  |u_\e(t)|^2_{L^2(D_\e)}+\e^2|\gamma_\e u
\e(t)|^2_{L^2(\p S_\e)}
 \big)+ \nonumber \\
&& \int_0^t \mathbf{E}\big(
|u_\e(t)|^2_{H_\e^1(D_\e)}+\e^2|\gamma_\e u_\e(t)|^2_{H^{1/2}(\p
S_\e)}
 \big)ds \leq
(1+\mathbf{E}|z_0|^2_{X^0_\e})C_T
\end{eqnarray}
and
\begin{eqnarray}\label{est4}
\mathbf{E}\big\{\sup_{t\in[0, T]}
|u_\e(t)|^2_{L^2(D_\e)}+\e^2|\gamma_\e u_\e(t)|^2_{L^2(\p S_\e)}
\big\}\leq (1+\mathbf{E}|z_0|^2_{X^0_\e})C_T.
\end{eqnarray}
\end{coro}

\vskip 0.6cm

 We recall a probability concept. Let $z$ be a random variable
taking values in a Banach space $\mathcal{S}$, namely, $z:\; \Om \to
\mathcal{S}$. Denote by $\mathcal{L}(z)$ the distribution (or law)
of $z$. In fact, $\mathcal{L}(z)$ is a Borel probability measure on
$\mathcal{S}$ defined as \cite{PZ92}

$$
\mathcal{L}(z)(A)=\mathbb{P} \{\omega: z(\omega) \in A \},
$$
for every event (i.e., a Borel set) $A$ in the Borel
$\sigma-$algebra $\mathcal{B}(\mathcal{S})$, which is the smallest
$\sigma-$algebra containing all open balls in $\mathcal{S}$.

\vskip 0.5cm

 As stated in \S   \ref{s1},  for the SPDE (\ref{e1}) we aim at deriving
 an
 effective equation in the sense of probability. A solution  $u_\e$ may be regarded as a
   random variable  taking values in  $L^2(0, T; L^2(D_\e))$.
  So for a solution $u_\e$ of
 (\ref{e1})-(\ref{e4}) defined on $[0, T]$,  we
focus on the behavior of distribution of $u_\e$ in $L^2(0, T;
L^2(D_\e))$ as $\e\rightarrow 0$. For this purpose, the tightness
\cite{Dudley} of
  distributions is necessary. Note   that the
function space changes with $\e$, which is a difficulty for
obtaining  the tightness of   distributions.   Thus we will treat
$\{\mathcal{L}(u_\e)\}_{\e>0}$ as a family of distributions on
$L^2(0, T; L^2(D))$ by extending $u_\e$ to the whole domain $D$.
 Recall that the distribution (or law ) of $\tilde{u}_\e$ is defined as:
$$
\mathcal{L}(\tilde{u}_\e)(A)=\mathbb{P} \{\omega:
\tilde{u}_\e(\cdot, \cdot, \omega) \in A \}
$$
for Borel set $A$ in $ L^2(0, T; L^2(D))$.\\

First we define the following spaces which will be used in our
approach. For Banach space $U$ and $p>1$, define $W^{1,p}(0 ,T;
U)$ as the space of functions $h\in L^p(0, T; U)$ such that
$$
 |h|^p_{W^{1, p}(0, T; U)}=|h|^p_{L^p(0, T; U)}
+\Big|\frac{dh}{dt}\Big|^p_{L^p(0, T; U)}<\infty.
$$
And for any $\alpha\in(0, 1)$, define $W^{\alpha, p}(0, T; U)$ as
the space of function $h\in L^p(0, T; U)$ such that
$$
 |h|^p_{W^{\alpha, p}(0, T; U)}=|h|^p_{L^p(0, T; U)}+
 \int_0^T\int_0^T\frac{|h(t)-h(s)|^p_U}{|t-s|^{1+\alpha
p}}dsdt<\infty.
$$
For $\rho\in(0,1)$, we denote by $C^\rho(0, T; U)$ the space of
functions  $h:[0,T]\rightarrow \mathcal{X}$ that are H\"older
continuous with exponent $\rho$.

\begin{theorem}\label{tight}
{\rm (\textbf{Tightness of distributions})} Assume that $z_0=(u_0,
v_0)$ be a $\big(\mathcal{F}_0,
\mathcal{B}(X^0_\e)\big)$-measurable random variable which is
independent of $W(t)$ with $\mathbf{E}|z_0|^2_{X^0_\e} <\infty $.
Then for any $T>0$,   $(\mathcal{L}(\tilde{u}_\e))_\e$, the
distributions of $(\tilde{u}_\e)_\e$, is tight in $L^2 (0, T;
L^2(D))\cap C(0, T;H^{-1}(D))$.
\end{theorem}
\begin{proof}
Denote the projection $(u, v)\rightarrow u$ by $P$. By the result
of Corollary \ref{est0},
\begin{equation}\label{1}
\mathbf{E}|u_\e|_{L^2(0,T;H^1_\e(D_\e))}^2\leq C_T.
\end{equation}
Write $z_\e(t)$ as
$$
z_\e(t)=z_\e(0)-\int_0^t\mathcal{A}_\e
z_\e(s)ds+\int_0^tF_\e(s,x)ds+\int_0^tG_\e(s, x)dW(s).
$$
Then by (\ref{ae}) and (\ref{weak}), when $(h, 0)\in X_\e^1$, we
have the following estimate, for some positive constant $C>0$
independent of $\e$
\begin{eqnarray*}
&&\Big|\Big( -P\int_0^t\mathcal{A}_\e
z_\e(s)ds+P\int_0^tF_\e(s,x)ds, h \Big)_{L^2(D_\e)}\Big|\\
&\leq&\Big|\int_0^t a(Pz_\e(s), h)ds\Big|+\Big|\int_0^t\big( f(s,x),
h\big
)_{L^2(D_\e)}ds\Big|\\
&\leq&C\Big(\int_0^t|u_\e(s)|_{H_\e^1(D_\e)}ds+\int_0^t|f(s)|_{L^2(D)}ds\Big)|h|_{H_0^1(D_\e)}.
\end{eqnarray*}
Thus we have
\begin{equation}\label{2}
\mathbf{E}\Big|-P\int_0^t\mathcal{A}_\e
z_\e(s)ds+P\int_0^tF_\e(s,x)ds \Big |^2_{W^{1,2}(0, T;
H^{-1}(D_\e))}\leq C_T.
\end{equation}
Let $M_\e(s,t)=\int_s^tG_\e(s,x)dW(s)$.   By   Lemma 7.2 of
\cite{PZ92} and the H\"older inequality, we have
\begin{eqnarray*}
\mathbf{E}|PM_\e(s,t)|^4_{L^2(D_\e)}&\leq&
\mathbf{E}|PM_\e(s,t)|^4_{L^2(D)}\leq c\mathbf{E}\big( \int_s^t
|g_1(\tau)|^2_{\mathcal{L}_2^{Q_1}}d\tau \big)^2\\&\leq&
K(t-s)\int_s^t\mathbf{E}|g_1(\tau)|^4_{\mathcal{L}_2^{Q_1}}d\tau\\
&\leq& K'(t-s)^2
 \end{eqnarray*}
for $t\in [s, T]$, and for positive constants $K$ and $K'$
independent of $\e$, $s$ and $t$. Therefore
\begin{equation}\label{3}
\mathbf{E}\int_0^T|PM_\e(0,t)|^4_{L^2(D_\e)}dt\leq C_T
\end{equation}
and for $\alpha\in (\frac{1}{4}, \frac{1}{2})$
\begin{eqnarray}\label{4}
\mathbf{E}\int_0^T\int_0^T\frac{|PM_\e(0,t)-PM_\e(0,s)|^4_{L^2(D_\e)}}{|t-s|^{1+4\alpha}}dsdt\leq
C_T.
\end{eqnarray}
Combining the estimates (\ref{1})-(\ref{4}) with the Chebyshev
inequality \cite{PZ92, Dudley}, it is clear that for any
$\delta>0$ there is a bounded set
$$
K_\delta\subset\mathcal{X}
$$
with $\mathcal{X}=L^2(0,T; H^1_\e(D))\cap \big(W^{1,2}(0,T;
H^{-1}(D))+W^{\alpha, 4}(0,T; L^2(D))\big)$, such that
\begin{equation*}
\mathbb{P}\{\tilde{u}_\e\in K_\delta \}>1-\delta.
\end{equation*}
Moreover by the compact embedding
$$
L^2(0,T; H^1(D))\cap W^{1,2}(0,T; H^{-1}(D))\subset L^2(0,T;
L^2(D))\cap C(0, T;H^{-1}(D))
$$
and
$$
L^2(0,T; H^1(D))\cap W^{\alpha, 4}(0,T; L^2(D))\subset L^2(0,T;
L^2(D))\cap C(0, T;H^{-1}(D)),
$$
we conclude that $K_\delta$ is compact in $L^2(0,T; L^2(D))\cap
C(0, T;H^{-1}(D))$. Thus $\{\mathcal{L}( \tilde{u}_\e)\}_\e$ is
tight in $L^2(0, T; L^2(D))\cap C(0, T;H^{-1}(D))$.

The proof is complete.
\end{proof}

\begin{remark}\label{fu}
When $f=f(t,x, u_\e)$ is nonlinear (i.e., it depends on $u_\e$)
but is also globally Lipschitz in $u_\e$, the results in Theorem
\ref{wellpose} and Corollary \ref{est0} still hold. For example,
see \cite{ChSch03} for such SPDEs with stochastic dynamical
boundary conditions.  Moreover, by the Lipschitz property, we have
$|f(t,x, u_\e)|_{L^2(D)}\leq C_T$. Hence a similar analysis as in
the proof of Theorem \ref{tight} yields the tightness of the
distribution for $u_\e$ in this globally Lipschitz nonlinear case.
This fact will be used in the beginning of \S \ref{s6} to get the
homogenized effective model when $f=f(t,x, u_\e)$ is globally
Lipschitz nonlinear.

In fact,  in \S \ref{s6}, we will also derive  homogenized
effective models   for three types of nonlinearities $f=f(t,x,
u_\e)$ that are \emph{not} globally Lipschitz in $u_\e$.
\end{remark}


\section{Two-scale convergence and some preliminary results}\label{s4}

In this section we present some basic results about the two-scale
convergence  \cite{All92, CD99}.

In the following we denote by $C^\infty_{per}(Y)$ the space of
infinitely differentiable functions in $\R^n$ that are periodic in
$Y$. We also   denote by $L^2_{per}(Y)$ or $H^1_{per}(Y) $ the
completion of $C^\infty_{per}(Y)$, in the usual norm of $L^2(Y)$
or $H^1(Y)$, respectively. We also introduce the space
$H^1_{per}(Y)/\R $, which is the space of the equivalent classes
of $u\in H^1_{per}(Y)$ under the following equivalent relation
$$
u \sim v \Leftrightarrow u-v=constant.
$$

\begin{defn}
A sequence of functions $u_\e(t,x)$ in $L^2(D_T)$ is said to be
two-scale convergent to a limit  $u(t,x,y)\in L^2(D_T\times Y)$,
if for any function $\varphi(x,y)\in C_0^\infty(D_T,
C^\infty_{per}(Y))$,
$$
\lim_{\e\rightarrow
0}\int_{D_T}u_\e(t,x)\varphi(t,x,\frac{x}{\e})dxdt
=\frac{1}{|Y|}\int_{D_T}\int_Yu_0(t,x,y)\varphi(t,x,y)dydxdt.
$$
This two-scale convergence is written as $u_\e
\overset{2-s}{\longrightarrow} u$.
\end{defn}

The following result ensures the existence of two-scale limit and
for the proof see \cite{All92,CD99}.
\begin{lemma}
Let $u_\e$ be a bounded sequence in $L^2(D_T)$. Then there exist a
function $u\in L^2(D_T\times Y)$  and a subsequence $u_{\e_k}$ with
$\e_k\rightarrow 0$ as $k\rightarrow \infty$ such that $u_{\e_k}$
two-scale converges to $u$.
\end{lemma}

\begin{remark}\label{wlimit}
Taking $\varphi$ independent of $y$ in the definition of two-scale
convergence, then $u_\e\overset{2-s}{\longrightarrow} u$ implies
that $u_\e$ weakly converges to its spatial average:
$$
u_\e(t,x) \rightharpoonup
\bar{u}(t,x)=\frac{1}{|Y|}\int_Yu(t,x,y)dy.
$$
So, we see that, for a given bounded sequence $L^2(D_T)$, the
two-scale limit $u(t,x,y)$ contains more information than the weak
limit $u(t,x)$: $u$ gives some knowledge on the periodic
oscillations of $u_\e$, while $\bar{u}$ is just the average with
respect to $y$. Another  advantage of the usage of two-scale
convergence is that we do not need an extending operator such as
in \cite{CD89, CDMZ91} in the homogenization procedure. For more
properties of two-scale convergence we refer to \cite{All92}.
\end{remark}

The following result is useful when considering   two-scale
convergence of the product of two convergent sequences, see
\cite{All92, CD99}.

\begin{lemma}\label{m2s}
Let $v_\e$ be a sequence in $L^2(D_T)$ that two-scale converges to
a limit $v(x,y)\in L^2(D_T\times Y)$. Further assume that

\begin{equation}\label{*}
\lim_{\e\rightarrow
0}\int_{D_T}|v_\e(t,x)|^2dxdt=\frac{1}{|Y|}\int_{D_T}\int_Y
|v(t,x,y)|^2dydxdt.
\end{equation}
Then, for any sequence $u_\e\in L^2(D_T)$, which two-scale
converges to a limit $u\in L^2(D_T\times Y)$, we have the weak
convergence of the product $u_\e v_\e$:
$$
u_\e v_\e \rightharpoonup \frac{1}{|Y|}\int_Yu(\cdot, \cdot ,
y)v(\cdot ,\cdot, y)dy,\;\; as \;\;\e\rightarrow 0 \;\; in \;\;
L^2(D_T).
$$
\end{lemma}

\begin{remark}
Condition (\ref{*}) always holds for a sequence of functions
$\varphi(t,x,x/\e)$, with $\varphi(t,x,y)\in L^2(D_T;
C_{per}(Y))$. Such functions $v_\e$ are called   admissible test
functions. With the additional condition (\ref{*}), the two-scale
convergence of $v_\e$ is also called strong two-scale convergence
\cite{All92}.
\end{remark}

Let $u_\e$ be a sequence of functions defined on $[0,T]\times
D_\e$ which is  bounded in $L^2(0,T; H_\e^1(D_\e))$. Then we have
the following result concerning the two-scale limit of the bounded
sequences $\tilde{u}_\e$ and $\widetilde{\nabla_x u}_\e$; for the
proof see \cite{All92}.

\begin{lemma}\label{2slimit}
There exist $u(t, x)\in H_0^1(D_T)$, $u_1(t,x,y)\in
L^2(D_T;H^1_{per}(Y))$ and a subsequence $u_{\e_k}$ with
$\e_k\rightarrow 0$ as $k\rightarrow \infty$,  such that
$$
\tilde{u}_{\e_k}(t,x)\overset{2-s}{\longrightarrow}\chi(y)u(t,x),
\;\; k\rightarrow\infty
$$
and
$$
\widetilde{\nabla_x u}_{\e_k}
\overset{2-s}{\longrightarrow}\chi(y)[\nabla_xu(t,x)+\nabla_yu_1(t,x,y)],\;\;\;\;
k\rightarrow\infty
$$
where $\chi(y)$ is the indicator   function  of $Y^*$ (which takes
value 1 on $Y^*$ and value 0 on $Y\setminus Y^*$).
\end{lemma}

\vskip 1cm

Since we consider the dynamical boundary condition, the technique
of transforming the surface integrals into the volume integrals is
useful in our approach. For this we follow the method of
\cite{Van81} (see also \cite{CD96}) for the nonhomogeneous Neumann
boundary problem for an elliptic equation.

For $h\in H^{-1/2}(\p S)$ and $Y$-periodic, define
$$
\Lambda_h=\frac{1}{|Y^*|}\int_{\p S}h(x)dx.
$$
 Also define
$$
\lambda_h=\frac{1}{|Y|}\langle h, 1 \rangle_{H^{{-1/2}},
H^{1/2}}=\vartheta \Lambda_h.
$$
Thus, in particular $\Lambda_1=\frac{|\p S|}{|Y^*|}$ and
\begin{equation}\label{coeff}
\lambda:=\lambda_1=\frac{|\p S|}{|Y|},
\end{equation}
where $|\cdot|$ denotes Lebesgue measure.

For  $h\in L^2(\p S)$ and $Y$-periodic, define $\lambda^\e_h\in
H^{-1}(D)$ as
$$
\langle \lambda^\e_h, \varphi\rangle=\e\int_{\p
S_\e}h(\frac{x}{\e})\varphi(x)dx,\;\;{\rm for}\;\; \varphi\in
H_0^1(D).
$$
Then we have the following result about the convergence of the
integral on the boundary.
\begin{lemma}\label{mue}
Let $\varphi_\e$ be a sequence in $H_0^1(D)$ such that
$\varphi_\e\rightharpoonup\varphi$ weakly in $H^1_0(D)$ as
$\e\rightarrow 0$. Then
$$
\langle \lambda_h^\e, \varphi_\e|_{D_\e} \rangle\rightarrow
\lambda_h\int_D\varphi dx, \;\;as\;\;\e\rightarrow 0.
$$
\end{lemma}

For the proof we refer to \cite{CD96}.


\section{Homogenized macroscopic model}\label{s5}

In this section we derive the effective macroscopic model for the
original model (\ref{e1}), by the two-scale convergence approach.
We first obtain a two-scale limiting model.  Then the homogenized
macroscopic model  is obtained by exploiting the relation between
weak limit and the two-scale limit.

By the proof of Theorem \ref{tight} for any $\delta>0$ there is a
bounded closed set $K_\delta\subset\mathcal{X}$ which is compact in
$L^2(0, T; L^2(D))$ such that
\begin{equation*}
\mathbb{P}\{\tilde{u}_\e\in K_\delta \}>1-\delta.
\end{equation*}
Then the Prohorov theorem and the Skorohod embedding theorem
(\cite{PZ92}) assure  that for any sequence $\{\e_ j\}_j$ with
$\e_j\rightarrow 0$ as $j\rightarrow \infty$, there exist
subsequence $\{\e_ {j(k)}\}$, random variables
$\{u^*_{\e_{j(k)}}\}\subset L^2(0, T; H_{\e_{j(k)}})$ and $u^*\in
L^2(0, T; H)$  defined on a new probability space $(\Omega^*,
\mathcal{F}^*, \mathbb{P}^*)$, such that
$$
\mathcal{L}(\tilde{u}^*_{\e_
{j(k)}})=\mathcal{L}(\tilde{u}_{\e_{j(k)}})
$$
and
$$
\tilde{u}^*_{\e_ {j(k)}}\rightarrow u^*\;\;in\;\;L^2(0,T; H)\;\; as
\;\;k\rightarrow \infty,
$$
for almost all $\omega\in\widehat{\Omega}$. Moreover $u^*_{\e_
{j(k)}}$ solves system (\ref{e1})-(\ref{e4}) with $W$ replaced by
Wiener process $W^*_k$ defined on probability space $(\Omega^*,
\mathcal{F}^*, \mathbb{P}^*)$ with same distribution as $W$.   In
the following, we will determine the limiting equation (homogenized
effective equation) that $u^*$ satisfies and the limiting equation
is independent of $\e$. After this is done we see that
  $\mathcal{L}(\tilde{u}_\e)$ weakly converges to $\mathcal{L}(u^*)$
  as $\e\downarrow 0$.\\

For $u_\e$ in set $K_\delta$, by Lemma \ref{2slimit} there is
$u(t,x)\in H_0^1(D_T)$ and $u_1(t,x,y)\in L^2(D_T; H^1_{per}(Y))$
such that

$$
\tilde{u}_{\e_ j}(t,x)\overset{2-s}{\longrightarrow}\chi(y)u(t,x)
$$
and
$$
\widetilde{\nabla_x u}_{\e_ j}
\overset{2-s}{\longrightarrow}\chi(y)[\nabla_xu(t,x)+\nabla_yu_1(t,x,y)].
$$
Then by Remark \ref{wlimit}
$$
\tilde{u}_{\e_
j}(t,x)\rightharpoonup\frac{1}{|Y|}\int_Y\chi(y)u(t,x)dy=\vartheta
u(t,x),\;\;{\rm weakly\;\; in}\;\; L^2(D_T).
$$
In fact by the compactness of $K_\delta$, the above convergence is
strong in $L^2(D_T)$. In the following, we will determine the
limiting equation, which is a two-scale system  that $u$ and $u_1$
satisfy. Then the limiting equation (homogenized effective
equation) that $u_0$ satisfies can be easily obtained by the
relation between weak limit and the two-scale limit.

Define a new probability space $(\Omega_\delta,\mathcal{F}_\delta,
\mathbb{P}_\delta)$ as
\begin{equation*}
\Omega_\delta=\{\omega\in\Omega: \tilde{u}_\e(\omega)\in K_\delta\},
\end{equation*}
\begin{equation*}
\mathcal{F}_\delta=\{F\cap\Omega_\delta: F\in\mathcal{F}\}
\end{equation*}
and
$$
\mathbb{P}_\delta(F)=\frac{\mathbb{P}(F\cap\Omega_\delta)}{\mathbb{P}(\Omega_\delta)},
\;\;{\rm for} \;\; F\in\mathcal{F}_\delta.
$$
Denote by $\mathbf{E}_\delta$ the expectation operator with respect
to $\mathbb{P}_\delta$. Now we restrict the system on the
probability space $(\Omega_\delta, \mathcal{F}_\delta,
\mathbb{P}_\delta)$.

Replace the test function $\varphi$ in (\ref{variational}) by
$\varphi_\e(t,x)=\phi(t,x)+\e \Phi(t,x,x/\e)$ with $\phi(t,x)\in
C_0^\infty(D_T)$ and  $\Phi(t,x,y)\in C_0^\infty(D_T;
C^\infty_{per}(Y))$. We will consider the terms in
(\ref{variational}) respectively .

By the choice of $\varphi_\e$ and noticing that
$\chi_{D_\e}\rightharpoonup \vartheta$, weakly$^*$ in $L^\infty(D)$,
we have
\begin{eqnarray}\label{t1}
&&\int_0^T\int_{D_\e}f(t,x)\varphi_\e(t,x) dxdt=\int_0^T\int_D
\chi_{D_\e}f(t,x)\varphi_\e(t,x) dxdt\nonumber \\
&\rightarrow& \vartheta\int_0^T\int_Df(t,x)\phi(t,x)dxdt,
\;\;\e\rightarrow 0.
\end{eqnarray}
And by the condition (\ref{gi})
\begin{eqnarray}\label{t2}
&&\int_0^T\int_{D_\e}g_1(t,x)\varphi_\e(t,x)dxdW_1(t)=\int_0^T\int_D\chi_{D_\e}g_1(t,x)\varphi_\e(t,x)dxdW_1(t)\nonumber\\
&\rightarrow& \vartheta\int_0^T\int_D g_1(t,x)\phi(t,x)dxdW_1(t),
\;\;\e\rightarrow 0\;\; in \;\; L^2(\Omega).
\end{eqnarray}
Integrating by parts and noticing that $\tilde{u}_\e$ converges
strongly to $\vartheta u(t,x)$ in $L^2(D_T)$,
\begin{eqnarray}
&&\int_0^T\int_{D_\e}\dot{u}_\e(t,x)\varphi_\e(t,x)
dxdt=-\int_0^T\int_{D_\e} u_\e(t,x)\dot{\varphi}_\e(t,x) dxdt\nonumber\\
&=&-\int_0^T\int_D\tilde{u}_\e(t,x)\dot{\phi}(t,x) dxdt-\e
\int_0^T\int_D \tilde{u}_\e(t,x)\dot{\Phi}(t,x, \frac{x}{\e}) dxdt\nonumber\\
&\rightarrow& -\int_0^T\int_D \vartheta u(t,x)\dot{\phi}(t,x)
dxdt=\int_0^T\int_D \vartheta \dot{u}(t,x)\phi(t,x) dxdt.\label{t3}
\end{eqnarray}
By the choice of $\varphi_\e$,
$$
\nabla_x\phi(t,x)+\nabla_y \Phi(t,x,
\frac{x}{\e})\overset{2-s}{\longrightarrow}
\nabla_x\phi(t,x)+\nabla_y \Phi(t,x,y), \;\; \e\rightarrow 0
$$
and
$$
\lim_{\e\rightarrow 0}\|\nabla_x\phi(t,x)+\nabla_y \Phi(t,x,
\frac{x}{\e})\|_{[L^2(D_T)]^n}=\frac{1}{|Y|}\int_{D_T\times
Y}\Big|\nabla_x\phi(t,x)+\nabla_y \Phi(t,x,y)\Big|^2dydxdt.
$$
Hence by  Theorem \ref{m2s}, we have
\begin{eqnarray}
&&\int_0^T\hspace{-0.2cm}\int_{D_\e}\hspace{-0.2cm}\nabla u_\e(t,x)
\nabla \varphi_\e(t,x)
dxdt=\int_0^T\hspace{-0.2cm}\int_{D_\e}\hspace{-0.2cm}\nabla
u_\e(t,x) (\nabla_x\phi(t,x)+\nabla_y \Phi(t,x, \frac{x}{\e})) dxdt\nonumber\\
&=&\hspace{-0.2cm}\int_0^T\hspace{-0.2cm}\int_D\hspace{-0.2cm}\widetilde{\nabla
u}_\e(t,x)(\nabla_x\phi+\nabla_y \Phi(t,x, \frac{x}{\e}))dxdt\nonumber\\
&\rightarrow
&\hspace{-0.2cm}\frac{1}{|Y|}\int_0^T\hspace{-0.2cm}\int_D\hspace{-0.05cm}\int_Y\chi(y)
\big[\nabla_xu(x)+\nabla_yu_1(x,y)\big]\big[\nabla_x\phi(t,x)+\nabla_y\Phi(t,x,y)
\big]dydxdt\nonumber \\
&=&\hspace{-0.2cm}\frac{1}{|Y|}\int_0^T\hspace{-0.2cm}\int_D\hspace{-0.05cm}\int_{Y^*}
\hspace{-0.2cm}\big[\nabla_xu(x)+\nabla_yu_1(x,y)\big]\hspace{-0.1cm}\big[\nabla_x\phi(t,x)+
\hspace{-0.1cm}\nabla_y\Phi(t,x,y)\big]dydxdt.\label{t4}
\end{eqnarray}
Now we consider the integrals on the  boundary. First for a fixed
$T>0$, it is easy to see that
\begin{eqnarray}
&&\e^2\int_0^T\int_{\p S_\e}\dot{u}_\e(t,x)\varphi_\e(t,x) dxdt\nonumber \\
&=&-\e^2\int_{\p S_\e}\int_0^T u_\e(t,x)\dot{\varphi}_\e(t,x)
dtdx\nonumber \\&=& -\e \Big\langle \lambda^\e_1, \int_0^T
\tilde{u}_\e(t,x)\dot{\varphi}_\e(t,x) dt\Big|_{D_\e} \Big \rangle
\rightarrow 0,\;\; \e\rightarrow 0. \label{t5}
\end{eqnarray}
And then
\begin{eqnarray}
&&\e b\int_0^T\int_{\p S_\e} u_\e (t,x)\varphi_\e(t,x) dxdt\nonumber \\
&=&\Big\langle \lambda^\e_1, \int_0^T \tilde{u}_\e(t,x)
\varphi_\e(t,x) dt\Big|_{D_\e} \Big
\rangle\nonumber \\
&\rightarrow& b\vartheta \lambda\int_0^T\int_Du(t,x)\phi(t,x)dxdt,
\;\;\e\rightarrow 0. \label{t6}
\end{eqnarray}
By the same method as above and the condition (\ref{gi}) we have
the limit of the stochastic integral on the boundary
\begin{eqnarray}
&&\e\int_0^T\int_{\p S_\e}g_2(t,x)\varphi_\e(t,x)dxdW_2(t)\nonumber\\
&\rightarrow&
\lambda\int_0^T\int_Dg_2(t,x)\phi(t,x)dxdW_2(t),\;\;\e\rightarrow
0,\;\;in\;\; L^2(\Omega). \label{t7}
\end{eqnarray}
Combining the above analysis in (\ref{t1})-(\ref{t7}) and by the
density argument we have
\begin{eqnarray}
&&\vartheta\int_0^T\int_D\dot{u}(t,x)\phi(t,x)dxdt\nonumber\\
&=&-\frac{1}{|Y|}\int_0^T\int_D\int_{Y^*}\big[\nabla_xu(x)+\nabla_yu_1(x,y)
\big] \big[ \nabla_x\phi(t,x)+\nabla_y\Phi(t,x,y)
 \big]dxdt\nonumber\\
&&-b\vartheta \lambda\int_0^T\int_Du(t,x)\phi(t,x)dxdt+\vartheta\int_0^T\int_Df(t,x)\phi(t,x)dxdt\nonumber\\
&&+\vartheta\int_0^T\hspace{-0.2cm}\int_D\hspace{-0.2cm}g_1(t,x)\phi(t,x)dxdW_1(t)+\lambda\int_0^T\hspace{-0.2cm}
\int_Dg_2(t,x)\phi(t,x)dxdW_2(t) \label{varl}
\end{eqnarray}
for any $\phi\in H_0^1(D_T)$ and $\Phi\in L^2(D_T;
H^1_{per}(Y)/\R)$. Integrating by parts, we see that (\ref{varl})
is the variational problem of the following two-scale homogenized
system
\begin{eqnarray}\label{2ssys}
\vartheta du=\big[-div_xA(\nabla_x u)-b\vartheta \lambda_1u+
\vartheta f\big]dt+ \vartheta g_1dW_1(t)+\lambda g_2dW_2(t),
\end{eqnarray}
\begin{equation}\label{2ssys1}
[\nabla_x u+\nabla u_1]\cdot\nu =0,\;\; on\; \p Y^*-\p Y
\end{equation}
where $\nu$ is the unit exterior norm vector on $\p Y^*-\p Y$ and
\begin{equation}
A(\nabla_x u)=\frac{1}{|Y|}\int_{Y^*} [\nabla_xu(t,x)+\nabla_y
u_1(t,x,y)]dy,
\end{equation}
with $u_1$ satisfying the following integral equation
\begin{equation}\label{u1}
\int_{Y^*}[\nabla_xu+\nabla_yu_1]\nabla_y\Phi
dy=0,\;\;u_1\;\;is\;\; Y-\mbox{periodic},
\end{equation}
for any $\Phi\in H_0^1(D_T;H_{per}^1(Y))$. The problem (\ref{u1})
has a unique solution for any fixed $u$, and so $A(\nabla_x u)$ is
well-defined. Furthermore $A(\nabla_xu)$ satisfies
\begin{equation}\label{A}
\langle A( \xi_1)-A( \xi_2), \xi_1-\xi_2 \rangle_{L^2(D),L^2(D)}\geq
\alpha\|\xi_1-\xi_2\|_{L^2(D)}^2
\end{equation}
and
\begin{equation}\label{AA}
|\langle A( \xi), \xi \rangle_{L^2(D),L^2(D)}|\leq
\beta\|\xi\|_{L^2(D)}^2
\end{equation}
 with some $\alpha$, $\beta>0$ and any $\xi$, $\xi_1$, $\xi_2\in H_0^1(D)$.
For   more detailed properties of $A(\nabla u)$ and (\ref{u1}) we
refer to \cite{FM87}. Then by the classical theory of the SPDEs,
  \cite{PZ92}, (\ref{2ssys})-(\ref{2ssys1}) is well-posed.

In fact $A(\nabla u)$ can be transformed to the classical
homogenized matrix by
\begin{equation}\label{u11}
u_1(t,x,y)=\sum_{i=1}^n\frac{\p u(t,x)}{\p
x_i}(w_i(y)-\mathbf{e}_i y),
\end{equation}
where $\{\mathbf{e}_i \}_{i=1}^n$ is the canonical basis of $\R^n$
and $w_i$ is the solution of the following cell problem (problem
defined on the spatial elementary cell)
\begin{eqnarray}
&&\D_y w_i(y)=0 \;\; in \;\; Y^*\\
&&w_i-\mathbf{e}_i y \;\;\;\;\mbox{is} \; Y-periodic\\
&&\frac{\p w_i}{\p \nu}=0 \;\; on\; \p S.
\end{eqnarray}
Then a simple   calculation yields
$$
A(\nabla u)=A^*\nabla u
$$
with $A^*=(A^*_{ij})$ being the classical homogenized matrix
 defined as
\begin{equation}\label{Aij}
A^*_{ij}=\frac{1}{|Y|}\int_{Y^*}w_i(y)w_j(y)dy.
\end{equation}
Then the above two-scale system (\ref{2ssys}) is equivalent  to
the following homogenized system,
\begin{eqnarray}
\vartheta du=\big[-div_x\big(A^*\nabla_x u\big)-b\vartheta \lambda
u+ \vartheta f\big]dt+ \vartheta g_1dW_1(t)+\lambda g_2dW_2(t).
\end{eqnarray}
Let $U(t,x)=\vartheta u(t,x)$. We thus have the limiting
homogenized equation
\begin{eqnarray}\label{limiting}
dU=\big[\hspace{-0.1cm}-\vartheta^{-1} div_x\big(A^*\nabla_x
U\big)\hspace{-0.1cm}-\hspace{-0.1cm}b \lambda U+ \vartheta
f\big]dt+ \vartheta g_1dW_1(t)+\lambda g_2dW_2(t).
\end{eqnarray}
And then $u^*$, we have mentioned at the beginning of this section,
satisfies (\ref{limiting}) with $W=(W_1, W_2)$ replaced by a Wiener
process $W^*$ with the same distribution as $W$. By the classical
existence result \cite{PZ92}, the homogenized model (\ref{limiting})
is well-posed. We formulate the main result of this section as
follows.

\begin{theorem} \label{macro}
 {\rm  (\textbf{Homogenized macroscopic model})} \\
Assume that (\ref{gi}) holds. Let $u_\e$ be the solution of
(\ref{e1})-(\ref{e4}). Then for any fixed $T>0$, the distribution
$\mathcal{L}(\tilde{u_\e})$ converges weakly to $\mu$ in $L^2(0, T;
H)$ as $\e\downarrow 0$, with $\mu$ being the distribution of $U$,
which is the solution of the following homogenized effective
equation
\begin{eqnarray}  \label{effect}
 dU=\hspace{-0.1cm}\big[\hspace{-0.1cm}-\vartheta^{-1} div_x\big(A^*\nabla_x U\big)-b \lambda U+
\vartheta f\big]dt+ \vartheta g_1dW_1(t)+\lambda g_2dW_2(t),
\end{eqnarray}
with the boundary condition $U=0$ on $\p D$, the initial condition
$U(0)=u_0/\vartheta$ and the effective matrix $A^*=(A^*_{ij})$
being determined by (\ref{Aij}). Moreover, the constant
coefficient $\vartheta=\frac{|Y^*|}{|Y|}$ is defined in the
beginning of \S 2 and $ \lambda=\frac{|\p S|}{|Y|}$ is defined in
(\ref{coeff}).
\end{theorem}

\begin{proof}
Noticing the arbitrariness of the choice of $\delta$, this is direct
result of the analysis of the first part in this section by the
Skorohod theorem and the $L^2(\Omega_\delta)-$convergence of $
\tilde{u}_\e$ on $(\Omega_\delta, \mathcal{F}_\delta,
\mathbb{P}_\delta)$.
\end{proof}

\vskip 0.5cm

\begin{remark} \label{rmk5.2}
It is interesting to note the following fact. Even when the
original microscopic model equation (\ref{e1}) is a deterministic
PDE (i.e., $g_1=0$), the homogenized macroscopic model
(\ref{effect}) is \emph{still} a stochastic PDE, due to the impact
of random dynamical interactions on the boundary of small scale
heterogeneities.
\end{remark}

\begin{remark}
For the macroscopic system (\ref{effect}),  we see that the fast
scale random fluctuations on the boundary is recognized or
quantified  in the homogenized equation, through the $
\mu_1g_2dW_2(t)$ term. The effect of random boundary evolution is
thus \emph{felt} by the homogenized system on the whole domain.
\end{remark}


\section{Homogenized macroscopic dynamics for nonlinear microscopic systems}\label{s6}

In this section, we derive homogenized macroscopic model for the
microscopic system (\ref{e1})-(\ref{e4}), when the nonlinearity
$f$ is either globally Lipschitz, or non-globally Lipschitz.

 As
Remark \ref{fu} has pointed out that if $f$ is a globally
Lipschitz nonlinear function of $u_\e$ all the estimates in $\S$
\ref{s3} hold. In fact, similar results in $\S$ \ref{s5} on
homogenized model also hold. In fact for    $f$ satisfying
$f(t,x,0)=0$ and
$$
|f(t,x,u_1)-f(t,x,u_2)|\leq L|u_1-u_2|
$$
for any $t\in\R$, $x\in D$ and $u_1$, $u_2\in\R$ with some
positive constant $L$. Since $\tilde{u}_\e\rightarrow  \vartheta
u$ strongly in $L^2(0, T; L^2(D))$ and by the Lipschitz property
of $f(t,x, u)$ with respect to $u$,  $f(t,x,
\tilde{u}_\e(t,x))\rightarrow f(t,x, u(t,x))$ strongly in $L^2(0,
T; L^2(D))$. (\ref{t1}) still hold for $f(t,x, u_\e)$. Then we can
obtain the same effective macroscopic system as (\ref{effect})
with nonlinearity $f=f(t, x, U)$:
\begin{eqnarray}  \label{effect888}
 dU=\hspace{-0.1cm}\big[\hspace{-0.1cm}-\vartheta^{-1} div_x\big(A^*\nabla_x U\big)-b \lambda U+
\vartheta f(t, x, U)\big]dt  \nonumber \\
+ \vartheta g_1dW_1(t)+\lambda g_2dW_2(t).
\end{eqnarray}

\medskip

For the rest of this section, we consider three types of nonlinear
systems with $f$   being non-global-Lipschitz   nonlinear function
in $u_\e$. The difficulty is at passing the limit $\e\rightarrow
0$ in the nonlinear term. These three types of nonlinearity
include: Polynomial nonlinearity; nonlinear term that is
sublinear; and nonlinearity that contains a gradient term $\nabla
u_\e$. We look at these nonlinearities case by case, and only
highlight the difference with the analysis in \S 5.

\bigskip
\textbf{Case 1: Polynomial nonlinearity }

 First we suppose $f$ is in the following
form
\begin{equation}\label{f}
f(t,x,u)=-a(t,x)|u|^pu
\end{equation}
with $0<a_0\leq a(t,x)\leq a_1$ for $t\in [0, \infty)$, $x\in D$.
And $p$ satisfies the following condition
\begin{equation}\label{p}
p\leq \frac{2}{n-2},\;\; {\rm if}\;\;n\geq 3;\;\;\; p\in \R,\;\;{\rm
if} \;\;n=2.
\end{equation}

For this case we need the following  \emph{Weak convergence lemma}
from Lions \cite{Li69}. \\

\emph{  Let $\mathcal{Q}$ be a bounded region in $\R\times \R^n$.
For any given functions $g_\e$ and $g$ in $L^p(\mathcal{Q})$
$(1<p<\infty)$, if
$$
|g_\e|_{L^p(\mathcal{Q})}\leq C,\;\; g_\e\rightarrow g\;\; {\rm in}
\;\;\mathcal{Q}\;\; {\rm almost\; everywhere}
$$
for some positive constant $C$, then $g_\e \rightharpoonup g$ weakly
in $L^p(\mathcal{Q})$.}

\bigskip

Noticing that  $F_\e(t, x, z_\e)=(f(t,x, u_\e), 0)$ and $ (F_\e(t,
x, z_\e), z_\e)_{X^0_\e}\leq 0, $  the results in Theorem
\ref{wellpose} can be obtained by the same method as in the proof
of Theorem \ref{wellpose}. Moreover by the assumption (\ref{p}),
$|f(t,x,u_\e)|_{L^2(D_T)}\leq C_T$, which by the analysis of
Theorem \ref{tight}, yields the tightness of the distribution of
$\tilde{u}_\e$.

Now we pass the limit $\e\rightarrow 0$ in $f(t,x,\tilde{u}_\e)$.
In fact, noticing that $\tilde{u}_\e$ converges strongly to
$\vartheta u$ in $L^2(0, T; L^2(D))$,   by the above weak
convergence lemma  with $g_\e=f(t,x,\tilde{u}_\e)$ and $p=2$,
$f(t,x,\tilde{u}_\e)$ converges weakly to $f(t,x, \vartheta u)$ in
$L^2(D_T)$. Threfore by the analysis for linear system in \S
\ref{s5}, we have the following result.

\begin{theorem}  \label{t6.1}
Assume that (\ref{gi}) holds. Let $u_\e$ be the solution of
(\ref{e1})-(\ref{e4}) with nonlinear term  $f$ being (\ref{f}).
Then for any fixed $T>0$,  the distribution
$\mathcal{L}(\tilde{u}_\e)$ converges weakly to $\mu$ in $L^2(0,
T; H)$ as $\e\downarrow 0$, with $\mu$ being the distribution of
$U$, which is the solution of the following homogenized effective
equation
\begin{eqnarray} \label{homo1}
 dU=\big[-\vartheta^{-1} div_x\big(A^*\nabla_x U\big)-b \lambda U+
\vartheta f(t,x, U)\big]dt   \nonumber \\
+ \vartheta g_1dW_1(t)+\lambda g_2dW_2(t),
\end{eqnarray}
with the boundary condition $U=0$ on $\p D$, the initial condition
$U(0)=u_0/\vartheta$ and the effective matrix $A^*=(A^*_{ij})$
being determined by (\ref{Aij}). Moreover, the constant
coefficient $\vartheta=\frac{|Y^*|}{|Y|}$ is defined in the
beginning of \S 2 and $ \lambda=\frac{|\p S|}{|Y|}$ is defined in
(\ref{coeff}).
\end{theorem}

\bigskip
\textbf{Case 2:  Nonlinear term that is sublinear }

More generally, we consider $f: [0, T]\times D\times\R \rightarrow
\R$ a measurable function which is continuous in $(x,\xi)\in
D\times\R$ for almost all $t\in [0,T]$ and satisfies
\begin{equation}\label{f1}
\big[f(t,x,\xi_1)-f(t,x,\xi_2)\big]\big[\xi_1-\xi_2\big]\geq
  0
\end{equation}
for $t\geq 0$, $x\in D$ and $\xi_1$, $\xi_2\in \R$. Moreover, we
assume that $f$ is sublinear,
\begin{equation}\label{f11}
|f(t,x,\xi)|\leq g(t)(1+|\xi|),\;\; \xi\in\R, \;\; t\geq 0,
\end{equation}
where $g\in L^\infty_{loc}[0, \infty)$. Notice that under the
assumption (\ref{f1}) and (\ref{f11}), $f$ may not be a Lipschitz
function.

 By the assumption (\ref{f11}) we can also have the tightness of the
 distributions of $\tilde{u}_\e$ and also conclude that
 $\chi_{D_\e}f(t,x,\tilde{u}_\e)$
 two-scale converges to a function denoted by $f_0(t,x,y)\in L^2(D_T\times
 Y)$. In the following we  need to identity $f_0(t,x, y)$.

 Let $\phi\in C_0^\infty(D_T)$ and  $\psi\in
C_0^\infty(D_T; C_{per}^\infty(Y))$. And for $\kappa>0$ let
\begin{equation}\label{xi}
\xi_\e(t,x)=\phi(t,x)+\kappa\psi(t,x,\frac{x}{\e}).
\end{equation}
Then by the assumption (\ref{f1}) we have
\begin{eqnarray*}
0&\leq& \int_0^T\int_{D_\e}\big[ f(t,x, u_\e)- f(t,x, \xi_\e)
\big]\big[u_\e-\xi_\e \big]dxdt \\ &=&
\int_{D_T}\chi\big(\frac{x}{\e}\big)\big[ f(t,x, \tilde{u}_\e)-
f(t,x, \xi_\e) \big]\big[\tilde{u}_\e-\xi_\e \big]dxdt\\
 & \overset{\triangle}= & I_\e=I_{1,\e}-I_{2,\e}-I_{3,\e}+I_{4,\e}
\end{eqnarray*}
with
\begin{eqnarray}\label{I1}
I_{1,\e}&=&\int_{D_T}\chi\big(\frac{x}{\e}\big)f(t,x,\tilde{u}_\e)\tilde{u}_\e dxdt\nonumber\\
&\overset{\e\rightarrow 0} \longrightarrow&
\frac{1}{|Y|}\int_{D_T}\int_Y f_0(t,x,y)\vartheta u(t,x) dydxdt,
\end{eqnarray}
\begin{eqnarray}\label{I2}
I_{2,\e}&=&\int_{D_T}\chi\big(\frac{x}{\e}\big)f(t,x,\tilde{u}_\e)\xi_\e dxdt\nonumber\\
&\overset{\e\rightarrow 0} \longrightarrow&
\frac{1}{|Y|}\int_{D_T}\int_Yf_0(t,x,y)\big[\phi(t,x)+\kappa\psi(t,x,y)\big]dydxdt,
\end{eqnarray}
\begin{eqnarray}\label{I3}
I_{3,\e}&=&\int_{D_T}\chi\big(\frac{x}{\e}\big)f(t,x,\xi_\e)\tilde{u}_\e dxdt\nonumber\\
&\overset{\e\rightarrow 0} \longrightarrow&
\frac{1}{|Y|}\int_{D_T}\int_Y\chi(y)f(t,x,\phi(t,x)+\kappa\psi(t,x,y))\nonumber
\\ && \vartheta u(t,x)dydxdt,
\end{eqnarray}
and
\begin{eqnarray}\label{I4}
I_{4,\e}&=&\int_{D_T}\chi\big(\frac{x}{\e}\big)f(t,x,\xi_\e)\xi_\e dxdt\nonumber\\
&\overset{\e\rightarrow 0} \longrightarrow&
\frac{1}{|Y|}\int_{D_T}\int_Y\chi(y)f(t,x,\phi(t,x)+\kappa\psi(t,x,y))\times \nonumber \\
&&\big[\phi(t,x)+\kappa\psi(t,x,y)\big]dydxdt.
\end{eqnarray}
In (\ref{I1})-(\ref{I4}) we have used the fact of strong two-scale
convergence of $\chi(\frac{x}{\e})$ and $f(t,x, \xi_\e)$, and the
strong convergence of $u_\e$ to $\vartheta u$.

Now we have
\begin{eqnarray*}
&&\lim_{\e\rightarrow 0}I_\e\\ &=&\int_{D_T}\int_Y \big[
f_0(t,x,y)-\chi(y)f(t,x,\phi+\lambda\psi)\big]\big[\vartheta
u(t,x)-\phi(t,x)-\kappa\psi \big]dydxdt\geq 0,
\end{eqnarray*}
for any $\phi\in C_0^\infty(D_T)$ and $\psi\in
C_0^\infty(D_T;C_{per}(Y))$. Letting  $\phi\rightarrow \vartheta u$,
dividing the above formula by $\kappa$  on both sides of the above
formula and letting $\kappa\rightarrow 0$ yields
$$
\int_{D_T}\int_Y \big[ f_0(t,x,y)-\chi(y)f(t,x,\vartheta u)\big]
\psi dydxdt\leq 0
$$
for any $\psi\in C_0^\infty(D_T; C_{per}(Y))$, which means
$$
f_0(t,x,y)=\chi(y)f(t,x,\vartheta u).
$$

Then by the similar analysis for linear systems in \S \ref{s5}, we
have the following homogenized model.

\begin{theorem} \label{t6.2}
Assume that (\ref{gi}) holds. Let $u_\e$ be the solution of
(\ref{e1})-(\ref{e4}) with nonlinear term f satisfying (\ref{f1})
and (\ref{f11}). Then for any fixed $T>0$, the distribution
$\mathcal{L}(\tilde{u}_\e)$ converges weakly to $\mu$ in $L^2(0,
T; H)$ as $\e\downarrow 0$, with $\mu$ being the distribution of
$U$, which is the solution of the following homogenized effective
equation
\begin{eqnarray} \label{Homo2}
 dU=\big[-\vartheta^{-1} div_x\big(A^*\nabla_x U\big)-b \lambda U+
\vartheta f(t,x, U)\big]dt  \nonumber  \\
+ \vartheta g_1dW_1(t)+\lambda g_2dW_2(t),
\end{eqnarray}
with the boundary condition $U=0$ on $\p D$, the initial condition
$U(0)=u_0/\vartheta$ and the effective matrix $A^*=(A^*_{ij})$
being determined by (\ref{Aij}).  Moreover, the constant
coefficient $\vartheta=\frac{|Y^*|}{|Y|}$ is defined in the
beginning of \S 2 and $ \lambda=\frac{|\p S|}{|Y|}$ is defined in
(\ref{coeff}).
\end{theorem}

\bigskip
\textbf{Case 3: Nonlinearity that contains a gradient term}

 We next consider $f$   in the following form containing a
gradient term,
\begin{equation}\label{f2}
f(t,x,u, \nabla u)=h(t,x,u)\cdot\nabla u
\end{equation}
where $h(t,x,u)=(h_1(t,x, u), \cdots, h_n(t,x,u))$ and each $h_i:
[0, T]\times D\times \R \rightarrow \R$, $i=1, \cdots, n$, is
continuous with respect to $u$ and
$h(\cdot,\cdot,u(\cdot,\cdot))\in L^2(0, T; L^2(D))$ for $u\in
L^2(0, T; L^2(D))$. Moreover assume that $h$ satisfies

\begin{enumerate}
  \item $|\langle h(t,x,u)\cdot\nabla u, v \rangle_{L^2}|\leq C_0 |\nabla
         u|_{L^2}|v|_{L^2}$ with some positive constant $C_0$.
  \item $|h_i(t,x,\xi_1)-h_i(t,x,\xi_2)|\leq k|\xi_1-\xi_2| $ for $\xi_1$, $\xi_2\in\R$, $i=1, \cdots, n$ and
  $k$ is a positive constant.
\end{enumerate}

 Now we have
\begin{equation}\label{fe2}
\big|(F_\e(t, x, z_\e), z_\e)_{X^0_\e}\big|=\big|(h(t,x, u_\e)\cdot
\nabla u_\e, u_\e)_{L^2}\big|\leq C_0|z_\e|_{X_\e^0}|z_\e|_{X_\e^1}.
\end{equation}
 By applying It\^o
formula to $|z_\e|_{X_\e^0}^2$, we obtain
\begin{eqnarray}\label{It2}
d|z_\e(t)|^2_{X^0_\e} +2(\mathcal{A}_\e z_\e,
z_\e)_{X^0_\e}dt&=&2(F_\e(t,x, z_\e),
z_\e)_{X^0_\e}dt+2(G_\e(t,x)dW(t), z_\e)_{X^0_\e}+\nonumber\\
&&|G_\e(t,x)|^2_{\mathcal{L}_2^Q}dt.
\end{eqnarray}
By (\ref{fe2}), coercivity (\ref{coercive}) of $a_\e(\cdot,\cdot)$
and the Cauchy inequality, integrating (\ref{It2}) with respect to
$t$ yields
\begin{eqnarray*}
&&|z_\e(t)|^2_{X^0_\e} +\bar{\alpha} \int_0^t|z_\e(s)|^2_{X^1_\e}ds\\
&\leq &|z_0|^2_{X^0_\e}+(2\bar{\beta}+\Lambda_1(\bar{\alpha}))\int_0^t|z_\e(s)|^2_{X_\e^0}ds+\\
&&2\int_0^t(G_\e(s) dW(s),
z_\e(s))_{X^0_\e}+\int_0^t|G_\e(s)|^2_{\mathcal{L}_2^Q}ds
\end{eqnarray*}
where $\Lambda_1$ is a positive constant depending on
$\bar{\alpha}$. Then by the Gronwall lemma we see that (\ref{est1})
and (\ref{est2}) hold. Moreover,   the fact
\begin{equation*}
|h(t,x, u)\cdot \nabla u|_{L^2}\leq C_0 |z_\e|_{X_\e^1},
\end{equation*}
together with the H\"older inequality yields
\begin{equation}
\mathbf{E}\Big|-P\int_0^t\mathcal{A}_\e z_\e(s)ds+P\int_0^tF_\e(s,x,
z_\e)ds \Big |^2_{W^{1,2}(0, T; H^{-1}(D_\e))}\leq C_T
\end{equation}
where $P$ is defined in Theorem \ref{tight}. Then by the same
discussion of Theorem \ref{tight}, we have the tightness of the
distributions of $\tilde{u}_\e$.

Now we   pass the limit $\e\rightarrow 0 $ in the nonlinear term
$f(t,x, u_\e, \nabla u_\e)$. In fact we restrict the system on
($\Omega_\delta, \mathcal{F}_\delta, \mathbb{P}_\delta$). By the
assumption (2) on $h$ and the fact that $\tilde{u}_\e$ strong
converges to $\vartheta u$ in $L^2(D_T)$, we have
\begin{equation*}
\lim_{\e\rightarrow 0}\int_{D_T}\Big[
h\big(t,x,\tilde{u}_\e(t,x)\big)-h\big(t,x, \vartheta u(t,x)\big)
\Big ]^2dxdt=0.
\end{equation*}

 For any $\psi\in C_0^\infty(D_T)$,
\begin{eqnarray}
&&\int_{D_T}h(t,x,\tilde{u}_\e)\cdot\widetilde{\nabla u_\e}\psi dxdt\nonumber \\
&=& \int_{D_T}\Big[h\big(t,x,\tilde{u}_\e\big)- h\big(t,x, \vartheta
u\big) \Big]\cdot\widetilde{\nabla
u_\e}\psi dxdt+\nonumber \\
&& \int_{D_T} h\big(t,x, \vartheta u\big)\cdot\widetilde{\nabla
u_\e}\psi dxdt\nonumber\\
&\overset{\e\rightarrow 0}\longrightarrow&
\frac{1}{|Y|}\int_{D_T}\int_Y h\big(t,x, \vartheta u\big)\cdot
\chi(y)\big[\nabla_xu+\nabla_y u_1 \big]\psi dydxdt.
 \label{nonlin2}
\end{eqnarray}
 Combining with the analysis for linear system in \S
\ref{s5}, we have the following result.

\begin{theorem} \label{t6.3}
Assume that (\ref{gi}) holds. Let $u_\e$ be the solution of
(\ref{e1})-(\ref{e4}) with nonlinear term (\ref{f2}). Then for any
fixed $T>0$, the distribution $\mathcal{L}(\tilde{u}_\e)$
converges weakly to $\mu$ in $L^2(0, T; H)$ as $\e\downarrow 0$,
with $\mu$ being the distribution of $U=\vartheta u$ which
satisfies  the following  homogenized effective equation
\begin{eqnarray}  \label{Homo3}
 dU=\big[- \vartheta^{-1} div_x\big(A^*\nabla_x U\big)-b\lambda U+
f^*(t,x,U,\nabla_x U)\big]dt   \nonumber \\
+ \vartheta g_1dW_1(t)+\lambda g_2dW_2(t),
\end{eqnarray}
where the boundary condition $U=0$ on $\p D$, the initial
condition $U(0)=u_0/\vartheta$, the effective matrix
$A^*=(A^*_{ij})$ is determined by (\ref{Aij}) and $f^*$ is the
following spatial average
$$
f^*(t,x,U,\nabla_x U)=\frac{1}{|Y|}\int_Y h\big(t,x, U\big)\cdot
\chi(y)\big[\vartheta^{-1}\nabla_x U+\nabla_y u_1 \big]dy
$$
with $u_1$ being given by (\ref{u11}) and $\chi(y)$ the indicator
function of $Y^*$. Moreover, the constant coefficient
$\vartheta=\frac{|Y^*|}{|Y|}$ is defined in the beginning of \S 2
and $ \lambda=\frac{|\p S|}{|Y|}$ is defined in (\ref{coeff}).
\end{theorem}

 \begin{remark}
 All the results in this paper hold when $\D$ is replaced by a more
 general strong elliptic operator {\rm div}$(A_\e\nabla u)$, where
 $A_\e$ is $Y-$periodic and satisfies the strong ellipticity
 condition.
 \end{remark}


\end{document}